\documentclass[12pt]{amsproc} 
\usepackage{amssymb,amsmath}
\usepackage{mathrsfs}

\newtheorem{definition}{Definition}[section]
\newtheorem{theorem}{Theorem}[section]
\newtheorem{lemma}{Lemma}[section]

\newtheorem{cor}{Corollary}[section]

\newtheorem{example}{Example}[section]

\def\Proof{{\bf Proof.} \begin{em} }
\def\endProof{\end{em}\hfill $\blacksquare$\par}

\def\bbq{\mathbb{Q}}
\def\bbr{\mathbb{R}}
\def\bbk{\mathbb{K}}
\def\bbl{\mathbb{L}}

\def\then{\rightarrow}
\def\iff{\leftrightarrow}

\def\scr#1{\mathscr{#1}}

\def\cb{\mathcal{B}}
\def\ci{\mathcal{I}}

\def\borel{{\mathcal{B}_+}}

\def\bigcups{\bigcup}

\def\<{\langle}
\def\>{\rangle}

\def\Ralowski{Ra{\l}owski}

\begin{document}
\title{Nonmeasurability in Banach spaces}
\author{Robert Ra{\l}owski}

\address{Robert \Ralowski, Institute of Mathematics, Wroc{\l}aw University of Technology, Wybrze\.ze Wyspia\'n\-skie\-go 27, 50-370 Wroc{\l}aw, Poland.}
\email{robert.ralowski@pwr.wroc.pl}


\begin{abstract}
We show that for a  $\sigma $-ideal $\ci$ with a Steinhaus property defined on Banach space, if two non-homeomorphic Banach with the same cardinality of the Hamel basis then there is a 
$\ci$ nonmeasurable subset as image by any isomorphism between of them. Our results generalize   results from \cite{kulka}.
\end{abstract}

\maketitle

\section{Notation and Terminology}

Throughout this paper, $X$, $Y$ will denote uncountable Polish spaces and $\cb(X)$ the Borel $\sigma $-algebra of $X.$  
We say that the ideal $\ci$ on $X$ has {\it Borel base}  if every  element $A\in\ci$ is contained in a Borel set in $\ci.$ (It is assumed that an ideal is always proper.)
The ideal consisting of all countable subsets of $X$ will be denoted by $[X]^{\le\omega }$ and the ideal of all meager subsets of $X$ will be denoted by $\bbk.$
Let $\mu $ be a continuous probability measure on $X.$ The ideal consisting of all $\mu $-null sets will be denoted by $\bbl_\mu.$ By the following well known result, $\bbl_\mu $ can be identified with the $\sigma $-ideal of Lebesgue null sets.

\begin{theorem}[\cite{SMS}, Theorem 3.4.23]
If $\mu $ is a continuous probability on $\cb (X),$ then there is a Borel isomorphism $h:X\rightarrow [0,1]$ such that for every Borel subset $B$ of $[0,1],$ $\lambda(B)=\mu(h^{-1}(B)),$ where $\lambda $ is a Lebesgue measure. 
\end{theorem}
 
\begin{definition}
We say that $(Z,\ci)$ is Polish ideal space if $Z$ is Polish uncountable space and $\ci$ is a $\sigma $-ideal on $Z$ having Borel base and containing all singletons. In this case, we set
$$\borel (Z)=\cb (Z)\setminus\ci.$$
A subset of $Z$ not in $\ci $ will be called a $\ci $-positive set; sets in $\ci $ will also be called $\ci $-null. Also, the $\sigma $-algebra generated by $\cb (Z)\cup\ci$ will be denoted by $\overline{\cb}(Z),$ called the $\ci $-completion of $\cb(Z).$
\end{definition}

It is easy to check that $A\in\overline{\cb }(Z)$ if and only if there is an $I\in\ci $ such that $A\bigtriangleup I$ (the symmetric difference) is Borel.

\begin{example}
Let $\mu $ be a continuous probability measure on $X.$ Then $(X,[X]^{\le\omega})$,
$(X,\bbk )$, $(X,\bbl_\mu )$ are Polish ideal spaces. 
\end{example}

\begin{definition}
A Polish ideal group is 3-tuple $(G,\ci,+)$ where $(G,\ci )$ is Polish ideal space and $(G,+)$ is an abelian topological group with respect to the Polish topology of $G.$
\end{definition}

Now we are ready to recall the crucial property which was introduced by Steinhauss see \cite{steinhaus}.
\begin{definition}\label{steinhaus} Let $(X,+)$ be any topological group with topology $\tau$. We say that ideal $\ci\subset\scr{P}(X)$ have Steinhaus property iff
$$
\forall A_1,A_1\in \scr{P}(X)\exists B_1,B_2\in\borel(X)\exists U\in\tau\;\;\; B_1\subset A_1\land B_2\subset A_2\land U\subset A_1+A_2.
$$
\end{definition}
In the same paper \cite{steinhaus} was proven that ideal of null sets poses the Steinhaus property.

\begin{definition} Let $(X,+)$ be any topological group and let $\ci\subset\scr{P}(X)$ be any invariant $\sigma$-ideal with singletons then $\ci$ has strong Steinhaus property iff
$$
\forall B\in\borel\forall A\in\scr{P}(X)\setminus\ci\;\; int(A+B)\ne\emptyset.
$$
\end{definition}

It is well known that ideal of meager sets in any topological group $(G,+)$ has strong Steinhaus property see \cite{shane}. Moreover the ideal of the null sets respect to Haar measure on the locally compact topological group $(G,+)$ has also strong Steinhaus property see \cite{BCS}.

\begin{definition}
Let $(X,\ci )$ be a Polish ideal space and $A\subseteq X$. We say that $A$ is $\ci$--nonmeasurable, if $A\notin\overline{\cb }(X).$ Further, we say that $A$ is completely $\ci $--nonmeasurable if
$$\forall B\in\borel (X) \;\;A\cap B\ne\emptyset\land A^c\cap B\ne\emptyset .
$$
\end{definition}

Clearly every completely $\ci $--nonmeasurable set is $\ci $--nonmeasurable. In the literature, completely $[X]^{\le\omega }$--nonmeasurable sets are called  Bernstein sets.
Also, note that  $A$ is completely $\bbl_\mu$--nonmeasurable if and only if the inner measure of $A$ is zero and the outer measure one.

For any set $E$, $|E|$ will denote the cardinality of $E.$

The rest of our notations and terminology are standard. For other notation and terminology in Descriptive Set Theory we follow \cite{SMS}.  

The main motivation of this paper is the Theorem about the nonmeasurability of the images under isomorphism over $\bbq$ between $\bbr^n$ and $\bbr^m$ whenever $m\ne n$. More precisely in \cite{kulka} the following Theorem was proved
\begin{theorem} Let $I$ is nontrivial invariant $\sigma$-ideal of subsets of the group $(\bbr,+)$ which has strong Steinhaus property such that $[\bbr]^\omega\subset I$. Let $f:\bbr^2\to\bbr$ be any linear isomorphism over $\bbq$. Then
\begin{enumerate}
 \item If $A\subset\bbr^2$ is a bounded subset of the real plane such that $A+T=\bbr^2$ for some $T\in [\bbr^2]^\omega$ then $f[A]$ is $I$-nonmeasurable subset of the real line.
 \item If $A\subset\bbr^2$ is a bounded subset of the real plane with non empty interior $int(A)\ne\emptyset$. Then the set $f[A]$ is completely $I$--nonmeasurable subset of the real line.
\end{enumerate}
\end{theorem}

\section{Results }

Here we present the main results of this paper.

\begin{theorem}\label{banach} Let $X$, $Y$ be a Banach spaces and let us suppose that
\begin{enumerate}
 \item $\ci\subset\scr{P}(Y)$ be an $\sigma$ - ideal with Steinhaus property,
 \item $\forall n\in\omega\forall A\in \ci\;\;n\ne 0 \then nA\in \ci$,
 \item $f:X\to Y$ by any isomorphism between $X$, $Y$ which is not homeomorphism.
\end{enumerate}
Then the image of the unit ball $f[K]$ is $\ci$ non-measurable in $Y$, where $K$ a unit ball of the space $X$ with the center equal $0\in X$.
\end{theorem}

\begin{theorem}\label{baire} Let $X,Y$ be any Banach spaces for which there exists linear isomorphism $f:X\to Y$ which is not continuous then image of the unit ball $K\subset X$ by $f$ has not Baire property.
\end{theorem}


\begin{theorem}\label{banach_kappa} Let $X$, $Y$ be a Banach spaces and let us assume that
\begin{enumerate}
 \item $\ci\subset\scr{P}(Y)$ be an $\kappa$ - complete additive invariant ideal with Steinhaus property,
 \item let $min\{ |S|:\;\; S\in\scr{P}(X)\land S \text{ is dense set in }X\}<\kappa$,
 \item $f:X\to Y$ by any isomorphism between $X$, $Y$, which is not homeomorphism.
\end{enumerate}
Then the image of the unit ball $f[K]$ is $\ci$ non-measurable in $Y$, where $K$ a unit ball of the space $X$ with the center equal $0\in X$.
\end{theorem}

\section{Proofs o the main results}

\Proof of the Theorem \ref{banach}. Let assume that image $f[K]$ is $\ci$-measurable in space $Y$. First of all let observe that
$$
Y=f[X]=f[\bigcups_{n=1}^\infty nK]=\bigcups_{n=1}^\infty f[nK]=\bigcups_{n=1}^\infty nf[K]
$$
then we have that $f[K]\not\in \ci$ is $\ci$ positive. Then by Steinhaus property of the ideal $\ci$ there exist unit ball $0\in U$ consisting $0$ for which $U\subset f[K]-f[K]=f[K-K]$ is hold. But $f$ is linear isomorphism then $f^{-1} [U]\subset K-K$ then $f^{-1}$ is bounded linear operator (so is continuous) then by Banach Theorem on invertible operator we have that $f$ is bounded then $f$ is homeomorphism what is impossible by $(3)$.
\endProof

Here we give the following Lemma:
\begin{lemma}\label{cathegory} Let $(X,\|\cdot\|)$ be any Banach space then $\sigma$-ideal $\bbk$ of the meager sets has Steinhaus property.
\end{lemma}
\Proof Let $A,B\in\scr{P}(X)\setminus\bbk$ be any meager positive subsets of the Banach space then there exists positive radius $r>0$ of two balls $K_1=K(x,\frac{r}{2}),K_2=K(y,r)\subset X$ with the centers $x,y\in X$ such that $A_1=A\cap K_1$ and $B_1=B\cap K_2$ are comeager subsets in the balls $K_1(x,r)$ and $K_2(y,r)$ respectively. Let $z=y-x$ and $r_0=\frac{r}{2}$ and let $K_0=K(z,r_0$ then we will show that
$$
\forall t\in K_0\;\; (t+A_!)\cap B_1\ne\emptyset
$$
First let observe that if $t\in K_0$ and $s\in K_1$ then
\begin{align*}
\|(t+s)-y\|=&\| (t-z+z)+(s-x+x)-y\|=\| (t-z)+(s-x)+z+x-y\|\\
=&\| (t-z)+(s-x)\|\le \| t-z\|+\| s-x\|<r_0+\frac{r}{2}=r
\end{align*}
thus we have $K_0+K_1\subset K_2$ and then $(K+A_1)\cap B_1$ is comeager in the open set $(K_0+K_1)\cap K_2$ then $(K+A_1)\cap B_1$ is nonempty set.

Thus there exists open set $K_0$ such that 
$$
\forall t\in K_0\;\; (t+A)\cap B\ne\emptyset
$$
Finally we have
$$
\forall t\in K_0\;\; (t+A)\cap B\ne\emptyset\iff \forall t\in K_0\exists a\in A\exists b\in B t+a=b\iff \forall t\in K_0 t\in B-A
$$
Thus we have $K_0\subset B-A$what finishes proof of this lemma.
\endProof

By the above Lemma we can give the Theorem \ref{baire}.

\Proof (of the Theorem \ref{baire}). By Lemma \ref{cathegory} the ideal of the meager sets of the space $Y$ has Steinhaus property and let observe that for any $\alpha\in\bbr\setminus\{ 0\}$ the map
$$
Y\ni y\mapsto \alpha\cdot y\in Y
$$
is the homeomorphism then second condition of the Theorem \ref{banach} is fullfiled then we are getting assertion.
\endProof

Immediately we have
\begin{cor} Let $X<\|\cdot\|)$ be infinite dimensional Banach space then there exists linear automorphism of $X$ for which image of the unit ball do not has Baire property.
\end{cor}
\Proof Let $\scr{B}$ is Hamel base which is subset of the unit ball of the space $X$. Let $\scr{B}_0=\{ e_n\in X:\; n\in\omega\}\subset\scr{B}$ be any countable subset of the our Hamel base $\scr{B}$. Let define $g:\scr{B}\to\scr{B}$ as follows:
$$
g(x)=\left\{\begin{matrix}
             (n+1)\cdot x	&x=e_n\in\scr{B}_0\\
		x		&x\in\scr{B}\setminus\scr{B}_0
            \end{matrix}\right. 
$$
Now we are ready to define $f:X\to X$ let $x\in X$ and let $A\in[\scr{B}]^{<\omega}$ and assume that
$$
x=\sum\limits_{e\in A}\alpha_e \cdot e
$$
then $f(x)=\sum_{e\in A}\alpha\cdot g(e)$. It is easy to see that $f$ is noncontinous linear automorphism of $X$. Then by Theorem \ref{baire} proof is finished.
\endProof

\Proof of the Theorem \ref{banach_kappa}. Let assume that image $f[K]$ is $\ci$-measurable in space $Y$. Let $\gamma<\kappa$ be such that $\gamma=|S|$ and the subset $S\subset X$ be dense in space $X$ then we have
$$
Y=f[X]=f[\bigcups_{x\in S} \{ x\}+K]=\bigcups_{x\in S} f[\{ x\}+K]=\bigcups_{x\in S} \{ f(x)\}+f[K]
$$
thus $f[K]\not\in \ci$ is $\ci$ positive. Then by Steinhaus property of the ideal $\ci$ there exist unit ball $0\in U$ consisting $0$ for which $U\subset f[K]-f[K]=f[K-K]$ is hold. But $f$ is linear isomorphism then $f^{-1} [U]\subset K-K$ then $f^{-1}$ is bounded linear operator (so is continuous) then by Banach Theorem on invertible operator we have that $f$ is bounded then $f$ is homeomorphism what is impossible by $(3)$.
\endProof

Immediately we are getting the following 
\begin{cor}\label{separable} Let $X$, $Y$ be a Banach spaces with the following properties:
\begin{enumerate}
 \item $X$ is separable Banach space,
 \item $\ci\subset\scr{P}(Y)$ be additive invariant $\sigma$-ideal with Steinhaus property,
 \item $f:X\to Y$ by any isomorphism between $X$, $Y$, which is not homeomorphism.
\end{enumerate}
Then the image of the unit ball $f[K]$ is $\ci$ non-measurable in $Y$, where $K$ a unit ball of the space $X$ with the centre equal $0\in X$.
\end{cor}

\noindent{\bf Acknowledgement} Author is very indebted to Professor Jacek Cicho{\'n} for help and critical remarks.


\vskip 1cm


\begin{thebibliography}{123}
\bibitem{BCS} A. Beck, H.H. Corson, A. B. Simon, The interior points of the product of two subsets of a locally compact group, Proc. Amer. Math. Soc. {\bf 9} (1958), 648--652
\bibitem{kulka} J. Cicho\'n, P. Szczepaniak, {\it When is the unit ball nonmeasurable?}, accepted to Fund. Math.
\bibitem{fremlin} D. H. Fremlin, Measure additive coverings and measurable selectors, Dissertationes Math. {\bf 260} (1987)
\bibitem{shane} E. J. McShane, Images of sets satisfying the condition of Baire, Ann. of Math. {\bf 51} (1950), 380--386
\bibitem{SMS} S.M. Srivastava, A course on Borel sets, Springer-Verlag, New York, 1998
\bibitem{steinhaus} H. Steinhaus, Sur les distances des points dans les ensembles de measure positive, Fund. Math. {\bf 1} (1920), 99--104
\end{thebibliography}
\end{document}